\documentclass[10pt]{amsart}

\usepackage{amsfonts}

\usepackage{graphicx}
\usepackage{amsmath}
\usepackage{setspace}

\newcommand{\halmos}{{\mbox{\, \vspace{3mm}}} \hfill
\mbox{$\Box$}}

\newtheorem{theorem}{Theorem}
\newtheorem{assumption}{Assumption}

\newtheorem{condition}[theorem]{Condition}
\newtheorem{conjecture}[theorem]{Conjecture}
\newtheorem{corollary}[theorem]{Corollary}

\newtheorem{definition}[theorem]{Definition}
\newtheorem{example}[theorem]{Example}

\newtheorem{lemma}[theorem]{Lemma}

\newtheorem{proposition}[theorem]{Proposition}
\newtheorem{remark}[theorem]{Remark}

\usepackage[english]{babel}
\usepackage[left=3cm,right=3cm,top=2.5cm,bottom=3cm]{geometry}

\usepackage{amsmath,amssymb,amsthm,bbm,color,graphics,version}
\usepackage{mathrsfs}

\newcommand{\bearno}{\begin{eqnarray*}}
\newcommand{\enarno}{\end{eqnarray*}}

\newcommand{\bs}{\boldsymbol}

\newcommand{\vb}{\vspace{3mm}}

\newcommand{\PP}{{\mathbb P}}
\newcommand{\EE}{{\mathbb E}}
\newcommand{\dd}{{\rm d}}
\newcommand{\vt}{\vartheta}
\newcommand{\s}{^\star}

\title{Quasi-stationary workload\\ in a L\'evy-driven storage system}

\author{Michel Mandjes}
\address{Korteweg-de Vries Institute for Mathematics,
University of Amsterdam, the Netherlands; E{\sc urandom},
Eindhoven University of Technology, the Netherlands; CWI,
Amsterdam, the Netherlands} \email{M.R.H.Mandjes@uva.nl}

\author{Zbigniew Palmowski}
\address{Mathematical Institute, University of Wroc{\l}aw, Poland.}
\email{zbigniew.palmowski@gmail.com}

\author{Tomasz Rolski}
\address{Mathematical institute, University of Wroc{\l}aw, Poland.}
\email{tomasz.rolski@gmail.com}

\thanks{This work is partially supported by the Ministry of Science and
Higher Education of Poland under the grants N N201 394137
(2009-2011).}

\date{\today}

\begin{document}

\begin{abstract}
In this paper we analyze the quasi-stationary workload of a
L\'evy-driven storage system. More precisely, assuming the system
is in stationarity, we study its behavior conditional on the event
that the busy period $T$ in which time 0 is contained has not
ended before  time $t$, as $t\to\infty$. We do so by first
identifying the double Laplace transform associated with the
workloads at time 0 and time $t$, on the event  $\{T>t\}.$ This
transform can be explicitly computed for the case of spectrally
one-sided jumps. Then asymptotic techniques for Laplace inversion
are relied upon to find the corresponding behavior in the limiting
regime that $t\to\infty.$ Several examples are treated; for
instance in the case of Brownian input, we conclude that the
workload distribution at time 0 and $t$ are both Erlang(2).

\vspace{3mm}


\end{abstract}

\maketitle

\pagestyle{myheadings} \markboth{\sc M.\ Mandjes --- Z.\ Palmowski
--- T.\ Rolski} {\sc Quasi-stationary workload in a L\'evy-driven storage
system}

\vspace{1.8cm}

\tableofcontents

\newpage

\section{Introduction}
Consider a storage system with L\'evy input, i.e.,  the process
$(Q(t))_t$ that evolves as a L\'{e}vy process $X(t)$ that is
reflected at $0$. In mathematical terms, this means that the
workload in the storage system at time $t$ is given by
\[Q(t):=X(t)-\inf_{s\leq t} X(s),\]
where $Q(0)=x$ for some initial workload $x\geq 0$. Assuming $\EE
X(1)<0$, there exists a stationary distribution $\pi$ of $Q(t)$;
it is seen that the stationary workload is distributed as the
all-time supremum:
\[\pi(x)=\PP\left(\sup_{t\geq 0} X(t)\leq x\right).\]
In the sequel we add the subscript $\pi$ to the probability
measure $\PP$ and the associated expectation $\EE$ when we wish to
indicate that $Q(0)$ is distributed according to this stationary
distribution.

Now let $T$ denote the {\it busy period}, that is $T=\inf\{t\ge0:\
Q(t)=0\}$. In this paper the object of our interest concerns the
existence and characterization of  the joint conditional
distribution
\[\label{def2}\lim_{t\rightarrow\infty}\PP_\pi(Q(0)\in \dd x, Q(t)\in \dd y\,|\,T >t)=:
\mu (\dd x, \dd y),\] where the convergence is to be understood in
the weak sense. We study this so-called {\it
quasi-stationary distribution} $\mu$ in detail; special attention
is paid to the marginal distributions $\mu(\cdot\times
\mathbb{R})=\mu^{\rm QS}_L(\cdot)$ and $\mu(\mathbb{R}\times
\cdot)=\mu^{\rm QS}_R (\cdot)$. As an aside we mention that
sometimes $\mu({\mathbb R}_+\times \dd y)$ is called
quasi-stationary distribution; here we do not follow that
convention.

\vb

A substantial body of work has been devoted to the analysis of
quasi-stationary distributions. Over the past decades, various
settings were considered; we here give a brief (non-exhaustive)
overview. Seneta and Vere-Jones \cite{Seneta}, Tweedie
\cite{tweedie}, Jacka and Roberts \cite{jackaroberts} focus on a
Markov chain setting, Iglehart \cite{Iglehart} addresses a random
walk setup, Kyprianou \cite{kyprianou} considers the M/G/1 queue
(i.e., a storage system with compound Poisson input), whereas
Martinez and San Martin \cite{Martmart} treat the case of Brownian
motion with drift. We also mention the contribution by Kyprianou
and Palmowski \cite{kyprpalm}, who found the quasi-stationary
distribution associated with a general light-tailed L\'{e}vy
process. Recently, Rivero \cite{rivero} (after appropriate
scaling) found the quasi-stationary distribution for the specific
situation that the L\'{e}vy process under study has a jump measure
with a regularly varying tail.

\vb

The contribution of this paper is twofold. In the first place, a
general formula for the double Laplace transform $(Q(0), Q(t))$ on
the event $\{T>t\}$, that is,
\begin{equation}\label{master}\int_0^\infty e^{-\vt t}
\EE_\pi[e^{-\alpha Q(0)-\beta Q(t)} , T>t]\,\dd t\end{equation}
is given (to which we refer to as the {\it master formula}).  The
derivation is based on the  Wiener-Hopf factorization, and can be
evaluated explicitly when all jumps are either all positive (the
so-called  spectrally-positive case) or all jumps are negative
(spectrally-negative case). These formulae allow us to identify
the quasi-stationary measures for the spectrally one-sided cases
(relying on Tauberian-type theorems), which can be regarded as the
second major contribution.

The paper provides interesting insights into the distribution of
the workload conditional on a long busy period. The distributions
found tend to be stochastically larger than the normal, stationary
distribution. For instance in the case of regulated standard
Brownian motion (with drift $-1$), both $\mu^{\rm QS}_L(\cdot)$
and $\mu^{\rm QS}_R(\cdot)$ correspond to Erlang(2) distributions
(with mean 2), whereas the stationary distribution is exponential
(with mean $\frac{1}{2})$. This type of insights can potentially
be used when setting up efficient importance sampling algorithms
\cite{GM} (as in those algorithms a change of measure is looked
for that mimics the distribution conditional on the rare event
under consideration).

\vb

This paper is organized as follows. In Section \ref{prel}
preliminaries are given: (i) we first recapitulate a set of main
results on fluctuation theory for L\'evy processes, and (ii) then
present Tauberian theorems that are useful in the context of this
paper (which can be used to identify the tail behavior of a random
variable from its Laplace transform). The main objective of
Section \ref{mf} concerns the derivation of the master formula,
i.e., an expression for (\ref{master}) in terms of the Wiener-Hopf
factorization (with explicit results for the spectrally one-sided
cases). Our findings on the quasi-stationary distribution are then
given in Section \ref{qsd}. The last section treats a number of
examples.

\section{Preliminaries}\label{prel}
\subsection{L\'evy processes}
Here we follow \cite{kyprianoubook} for definitions,  notations
and basic facts on L\'evy processes. Let in the sequel $X\equiv
(X(t))_t$ be a L\'{e}vy process which is defined on the filtered
space $(\Omega,\mathscr{ F},\{\mathscr{F}_t\}_{t\geq 0}, \PP)$
with the natural filtration that satisfies the usual assumptions
of right continuity and completion.  Later if we write $\PP_x$, it
means that $\PP_x(X(0)=x)=1$ and $\PP_0=\PP$; similarly, $\EE_x$
is expectation with respect to $\PP_x$. We denote by $\Pi(\cdot)$
jump measure of $X$. Later we will focus on asymmetric L\'{e}vy
processes, which are either spectrally negative (having
nonpositive jumps) or spectrally positive L\'{e}vy processes
(having nonnegative jumps).

\vb

{\it First passage times.} For any L\'{e}vy process we can define
its Laplace exponent $\psi({\eta})$ by
\begin{equation}\label{psi}
\EE e^{{\eta} X(t)}=e^{t\psi({\eta})},\end{equation} for
${\eta}\in \Theta$ such that the left hand side of
(\ref{psi}) is well-defined (from now on we will assume that that
this set $\Theta$ is not empty). Later, we also need the first
passage time \[\tau(x):=\min\{t: X(t)\ge x\},\] which for a
spectrally negative process $X$ with positive drift (i.e., $\EE
X(1)>0$) has Laplace transform
\begin{equation}\label{LTtau}\EE e^{-s\tau(x)}=e^{-\Phi(s)x}, \end{equation}
where $\Phi(s):=\sup\{{\eta}\ge0: \psi({\eta})=s\}$
is the right inverse of $\psi$  (see \cite{kyprianoubook} for details).

\vb

{\it Exponential change of measure.} For ${\eta} \in \Theta$
we define a new probability measure  $\PP_x^{{\eta}}$ by the
relation
\[\left. \frac{\dd \PP_{x}^{{\eta}}}{\dd \PP_{x}}\right| _{\mathscr{F}_{t}}=e^{{\eta}(X(t)-x)-\psi({\eta})t}\;;\]
we say that we have performed an exponential change of measure.
Under $\PP^{{\eta}}$, the process $X$ is still a L\'evy
process, but now with Laplace exponent
\begin{equation}\label{newpsi}\psi_{\eta}(\beta):=\psi({\eta}+\beta)-\psi({\eta}).\end{equation}
We will use the subscript ${\vartheta}$ to  indicate that the
quantity under consideration relates to $\PP^{\eta}$.

\vb

{\it Dual process.} We will also consider  the so-called {\it dual
process} $\hat{X}_t=-X_t$ with jump measure
$\hat{\Pi}\left(0,y\right)=\Pi\left(-y,0\right)$. Characteristics
of $\hat{X}$ will be indicated by using the same symbols as for
$X$, but with a `$\hat{\hspace{2mm}}$' added.

\vb

{\it Ladder heights.}
For the process $X$ we define the associated   $(L^{-1}(t),H(t))_{t}$:
\[L^{-1}(t):= \left \{ \begin{array}{ll}
\inf\{s>0: L(s)>t \} & \textrm{if $t < L(\infty)$,} \\ \infty & \textrm{otherwise,}
\end{array} \right.\]
and
\[H(t):= \left \{ \begin{array}{ll}
X_{L^{-1}(t)} & \textrm{if $t < L({\infty})$,} \\ \infty & \textrm{otherwise},
\end{array} \right.\]
where $L\equiv (L(t))_t$ is the local time at the  maximum
\cite[p.\ 140]{kyprianoubook}. Recall that $(L^{-1},H)$ is a
bivariate subordinator with the Laplace exponent
\[\kappa({\varphi},\beta):=-\frac{1}{t}\log \EE\left (e^{-{\varphi} L^{-1}(t)-\beta H(t)}\textbf{1}_{\{t\leq L({\infty})\}}\right)\] and with the jump measure
$\Pi_H$. In addition to this, we define the  {\it descending ladder  height
process} $(\hat{L}^{-1}(t),\hat{H}(t))_{t\geq 0}$ with the Laplace
exponent $\hat{\kappa}({\vartheta}, \beta)$ constructed from the
dual process $\hat{X}$. Recall that under the stability assumption
$\EE X(1)<0$, the random variable $L(\infty)$ has an exponential
distribution with parameter $\kappa(0,0)$. Moreover, for a
spectrally negative L\'{e}vy process  the Wiener-Hopf
factorization states that
\begin{equation}\label{kappy}
\kappa({\varphi}, \beta)=\Phi({\varphi})+\beta,\qquad \hat{\kappa}({\varphi}, \beta)=\frac{{\varphi}-\psi(\beta)}{\Phi({\varphi})-\beta};
\end{equation}
see \cite[p. 169-170]{kyprianoubook}. It follows that
$
\kappa(0,0)=\psi^\prime(0+).$

We introduce a {\it potential measure} $\mathscr{U}$ defined by
\[ \mathscr{U}(\dd x,\dd s)=\int_{t=0}^{\infty} \PP\left(L^{-1}(t) \in \dd s, H(t) \in \dd x\right) \dd t\]
with the Laplace transform
$ \int_{[0,\infty)^2} e^{-{\varphi} s -\beta x} \mathscr{U}(\dd x,\dd s)=1/\kappa({\varphi},\beta) $ and renewal function
\[V(\dd x)=\int_{s=0}^\infty\mathscr{U}(\dd x,\dd s)=\EE\left(\int_{t=0}^{\infty} \textbf{1}_{\{H(t)\in \dd x\}} \dd t \right).\]
In particular,
\begin{equation}\label{LTV}\int_0^\infty e^{-\beta x}V(x)\,\dd x =\frac{1}{\beta \kappa(0,\beta)}.
\end{equation}
For a spectrally negative L\'{e}vy process, the upward ladder
height process is just a linear drift, and hence the renewal
measure corresponds to the Lebesgue measure:
\begin{equation}\label{renfunleb}
V(\dd x)=\dd x.
\end{equation}
From (\ref{kappy}) we have \cite[p. 195]{kyprianoubook} that
\[
\int_0^\infty e^{-{\varphi} z}\hat{V}(\dd z)=\frac{{\varphi}}{\psi({\varphi})}.\]

\subsection{Tauberian-type results}\label{sec:sm}
\newcommand{\imi}{{\rm i}}

Consider a function $f:\mathbb{R}\to\mathbb{R}$ such that $f(z)=0$
for $\Re(z)<0$. Let $\tilde f(z):=\int_0^\infty e^{-z x}f(x)\,\dd
x$ be its Laplace transform. Consider its singularities; among
these, let ${\vartheta}{\s}<0$ the one  with the largest real
part. Notice that this yields the integrability of
$\int_0^\infty|f(x)|\,\dd x$. The inversion formula the reads
\[f(x)=\frac{1}{2\pi \imi}\int_{a-\imi \infty}^{a+\imi\infty}\tilde{f}(z) e^{zx}\,\dd z\]
for some (and then any) $a>{\vartheta}{\s}$.

\vb

We now focus on a class of theorems that infer the tail behavior
of a function from its Laplace transform, commonly referred to as
{\it Tauberian theorems}. Importantly, the behavior of the Laplace
transform around the singularity ${\vartheta}{\s}$ plays a crucial
role here. The following heuristic  principle given in \cite{AW1997} is
often relied upon. Suppose that for ${\vartheta}{\s}$, some
constants $K$ and $C$,  and a non-integer $s > 0$,
\[
\tilde{f}({\vartheta})=K - C
({\vartheta}-{\vartheta}{\s})^{s}+o(({\vartheta}-{\vartheta}{\s})^{s}),\qquad
\mbox{as ${\vartheta}\downarrow {\vartheta}{\s}$.}\] Then
\[
f(x)=\frac{C}{\Gamma(-s)} x^{-s-1}e^{{\vartheta}{\s}
x}(1+o(1)),\qquad \mbox{as $x\to\infty$,}\] where $\Gamma(s)$ is
the gamma function. Below we specify conditions under which this
relation can be rigorously proven. Later in our paper we apply it
for the specific case that $s=1/2$; recall that
$\Gamma(-1/2)=-2\sqrt{\pi}$.

\begin{comment}
%Unfortunately, as the following example shows, the utmost care is
%required; as mentioned above, conditions must be imposed for the
%above principle to hold. Consider the function $f(t)=\cos t$, with
%Laplace transform
%$\tilde{f}({\vartheta})={\vartheta}/(1+{\vartheta}^2)$. It has a
%singularity at ${\vartheta}{\s}=0$ and $\tilde f({\vartheta})\sim
%{\vartheta}$ as ${\vartheta}\to0$. Asymptotics of $f(t)$ (as
%$t\to\infty$), however, do not exist. Notice that
%$\tilde{f}({\vartheta})$ is analytic only in
%$\{z\in\mathbb{C}:\Re(z)>0\}$. \footnote{\tt MM: Is this really
%correct? Is there a pole at $\theta = 0$??}
\end{comment}
\vb

A formal justification of the above relation can be found in
Doetsch \cite[Theorem~37.1]{Doetsch1974}. Following Miyazawa and
Rolski \cite{miyazawaRolski}, we consider the following specific
form. For this we  first recall the concept of the $\mathfrak
W$-contour with an half-angle of opening $\pi/2<\psi\le \pi$, as
depicted on \cite[Fig. 30, p. 240]{Doetsch1974}; also, ${\mathscr
G}_{\zeta{\s}}(\psi)$  is the region between the contour
$\mathfrak W$ and the line $\Re(z)=0$. More precisely,
\[{\mathscr G}_{\alpha}(\delta) \equiv \{z \in \mathbb{C}; \Re(z) < 0, z \ne \alpha,
|\arg (z - \alpha)| < \delta \},\]  where $\arg z$ is the principal part of the argument of the complex number $z$.

In the following theorem, conditions  are identified such that the
above principle holds; we refer to this as the {\it Heaviside's
operational principle}, or simply {\it Heaviside principle}.

\begin{theorem}[Heaviside principle]\label{t.tauberian}
Suppose that for  $\tilde f:\mathbb{C}\to \mathbb{C}$ and
$\zeta{\s} < 0$ the following three conditions hold:
\begin{itemize}\item [(A1)] $\tilde f(\cdot)$ is analytic in a
region ${\mathscr G}_{\zeta{\s}}(\psi)$ for some $\pi/2<\psi\le
\pi$;
\item [(A2)] $\tilde f(z) \to 0$ as $|z| \to \infty$ for $z \in {\mathscr G}_{\zeta{\s}}(\psi)$;
\item [(A3)]  for some constants $K$ and $C$, and a non-integer $s>0$,
\begin{eqnarray}
\label{eqn:bridge case 1} \tilde f(z)=K\;{\boldsymbol 1}_{\{s >
0\}} - C (z-\zeta{\s})^{s}+o((z-\zeta{\s})^{s}),
\end{eqnarray}
where ${\mathscr G}_{\zeta{\s}}(\psi)\ni z \to \zeta{\s}$.
\end{itemize}
Then
\begin{eqnarray*}
f(x)=\frac{C}{\Gamma(-s)} x^{-s-1}e^{\zeta{\s} x}(1+o(1))\;,
\end{eqnarray*}
as $x\to\infty$, where $K:=\tilde f(\zeta{\s})$ if $s > 0$.
\end{theorem}

We now discuss when assumption (A1) is satisfied.  To check that
the Laplace transform $\tilde{f}(\cdot)$ is  analytic in the
region $\mathscr G_{\zeta{\s}}(\psi)$, we can use the concept of
semiexponentiality  of $f$  (see \cite{henrici}).
\begin{definition}[Semiexponentiality]
$f$  is said to be {\em semiexponential} if for some  $0< \phi\le\pi/2$, there exists finite
and strictly negative
$\gamma({\vartheta})$, defined  as the infimum of all such $a$ such that
\[\left|f(e^{\imi{\vartheta}}r)\right|<e^{ar}\]
for all sufficiently large $r$;
here $-\phi\le {\vartheta}\le \phi$
and $\sup \gamma({\vartheta})<0$.
\end{definition}

Relying on this concept, the following sufficient condition for (A1) applies.

\begin{proposition}\label{a1semexp}{\em \cite[Thm. 10.9f]{henrici}}
Suppose that $f$ is  semiexponential with $\gamma({\vartheta})$
fulfilling the following conditions: {\rm (i)}
 $\gamma=\gamma(0)<0$, {\em (ii)} $\gamma({\vartheta})\geq\gamma(0)$ in a neighborhood of ${\vartheta}=0$, and
{\rm (iii)} it is smooth. Then {\em (A1)} is satisfied.
\end{proposition}

\section{Master formula}\label{mf}

The objective of this section is to  derive a general formula for
the double Laplace-Stieltjes transforms
\[L_x({\vartheta};\alpha,\beta):=\int_0^\infty e^{-{\vartheta} t} \EE_x[e^{-\alpha x-\beta Q(t)},T>t]
\,\dd t\;\]
and
\[
L({\vartheta};\alpha,\beta):=\int_0^\infty
L_x({\vartheta};\alpha,\beta)\,\dd \PP(Q(0)\le x)= \int_0^\infty
e^{-{\vartheta} t} \EE_\pi[e^{-\alpha Q(0)-\beta Q(t)},T>t] \,\dd
t.\] Let $e_{{\vartheta}}$ be an exponentially distributed random
variable with parameter ${\vartheta}>0$, independent of the
process $X$. Denote
\[\overline{X}(t):=\sup_{s\leq t}X(s),\qquad \underline{X}(t):=\inf_{s\leq
t}X(s).\] Following the idea behind \cite[Th. VI.20]{bert96},  we
can prove the following result.

\begin{theorem}\label{lemma:basic}
For $\alpha,\beta,{\vartheta}>0$,
\[
L_x({\vartheta};\alpha,\beta)=\frac{1}{{\vartheta}}\frac{\kappa({\vartheta},0)}{\kappa({\vartheta},\beta)}
 e^{-(\alpha +\beta) x}\left(\int_0^x e^{\beta z}
\hat{\PP}(\overline {X}({e_{{\vartheta}}})\in \dd
z)\right).\]\end{theorem}

\proof It is elementary that
\begin{eqnarray*}
\lefteqn{\int_0^\infty e^{-{\vartheta} t} \EE_x[e^{-\alpha x-\beta
Q(t)},T>t] \,\dd t=\int_0^\infty e^{-{\vartheta} t}
\EE_x[e^{-\alpha x-\beta X(t)},T>t]
\,\dd t}\\
&=&\EE_x\int_0^Te^{-{\vartheta} t} e^{-\alpha x -\beta X(t)}\dd t= \frac{1}{{\vartheta}} \, \EE_x[ e^{-\alpha x-\beta X(e_{{\vartheta}})},e_{{\vartheta}}<T].
\end{eqnarray*}
Recall that \cite[Thm. 6.16(i)]{kyprianoubook}) under $\PP$ we have that $\underline{X}(e_{\vartheta})$ and $X(e_{\vartheta})-\underline{X}(e_{\vartheta})$
are independent; in addition $X(e_{\vartheta})-\underline{X}(e_{\vartheta}) =_{\rm d}\overline{X}(e_{\vartheta}).$
It thus follows that
\begin{eqnarray*}
\lefteqn{\frac{1}{{\vartheta}}  \,\EE_x [ e^{-\beta \underline{X}(e_{{\vartheta}})},e_{{\vartheta}}<T]=
\frac{1}{{\vartheta}} \, \EE_x[ e^{-\beta (X(e_{\vartheta})-\underline{X}(e_{{\vartheta}}))
-\beta
X(e_{\vartheta})},e_{{\vartheta}}<T]}\\
&=&\frac{1}{{\vartheta}}e^{-\beta x} \EE[e^{-\beta
(X(e_{\vartheta})-\underline{X}(e_{{\vartheta}}))}]
\,\EE[e^{-\beta
\underline{X}(e_{\vartheta})},\underline{X}(e_{\vartheta})>-x]\\
&=&\frac{1}{{\vartheta}}e^{-\beta x} \EE [e^{-\beta
\overline{X}(e_{{\vartheta}})}] \,\hat{\EE}[e^{\beta
\overline{X}(e_{\vartheta})}, \overline{X}(e_{\vartheta})<x]\;.
\end{eqnarray*}
Following \cite[Th. 6.16(ii)]{kyprianoubook}, we obtain that
\begin{equation}\label{supdual}\EE[e^{-\beta \overline{X}(e_{{\vartheta}})}]=\frac{\kappa({\vartheta},0)}{\kappa({\vartheta},\beta)}\end{equation}
and hence
$$\hat{\EE}[e^{\beta
\overline{X}(e_{\vartheta})};\overline{X}(e_{\vartheta})<x]=\int_0^xe^{\beta z}
\hat{\PP}(\overline{X}(e_{\vartheta})\in\, {\dd}z)$$
which completes the proof.
\halmos

\vb

We now evaluate $L({\vartheta};\alpha,\beta)$ for the spectrally one-sided cases.

\begin{proposition}\label{p.joint.transform} If $X$ is a spectrally positive L\'{e}vy process, then
\begin{equation}
L({\vartheta};\alpha,\beta)=\frac{\hat{\psi}'(0+)}{{\vartheta}-\hat{\psi}(\beta)}\left(\frac{\alpha+\beta}{\hat{\psi}(\alpha+\beta)}-\frac{\alpha+\hat{\Phi}({\vartheta})}
{\hat{\psi}(\alpha+\hat{\Phi}({\vartheta}))}\right).\label{eq.Lsp.pos}
\end{equation}
\end{proposition}

\proof
Note that
\[\hat{\PP}(\overline{X}(e_{\vartheta})\geq x)=\hat{\PP}(\tau(x)\le e_{\vartheta})=\hat{\EE}e^{-{\vartheta} \tau(x)}=e^{-\hat{\Phi}({\vartheta})x}.\]Integration by parts yields
\[
\int_{0}^xe^{\beta z}\,\hat{\PP}(\overline{X}(e_{\vartheta})\in
\dd z)
%&=&
%-e^{\beta z}\hat{P}_0(\overline{X}(e_{\vartheta})\ge z)\left|_{z=0}^{z=x}\right.+\beta
%\int_0^x e^{\beta z}\hat{P}_0(\overline{X}(e_{\vartheta})\ge z)\,dz\\
%&=&1-e^{\beta x}\hat{E}e^{-{\vartheta} \tau(x)}+\beta\int_0^x e^{\beta z}e^{-z\hat{\Phi}({\vartheta})}\,dz\\
%&=&1-e^{(\beta-\hat{\Phi}({\vartheta}))x}+\frac{\beta}{\beta-\hat{\psi}({\vartheta})}\left(e^{(\beta-\hat{\Phi}({\vartheta}))x}-1\right)\\
=\frac{\hat{\Phi}({\vartheta})}{\beta-\hat{\Phi}({\vartheta})}\left(e^{(\beta-\hat{\Phi}({\vartheta}))x}-
1\right)\;.\] The Pollaczek-Khintchine formula   \cite[Eqn.\
(4.14), p.\ 101]{kyprianoubook} states that
\begin{equation}\label{pollaczek}\tilde{\pi}(s):=\int_0^\infty  e^{-s x }\,\pi(\dd x)=\frac{\hat{\psi}'(0+) s}{\hat{\psi}(s)}.\end{equation}
It now follows that
\begin{eqnarray*}
\lefteqn{L({\vartheta};\alpha,\beta)=\frac{1}{{\vartheta}}\frac{\kappa({\vartheta},0)}{\kappa({\vartheta},\beta)}
\int_0^\infty e^{-(\alpha +\beta) x}\int_0^x e^{\beta z}
\hat{\PP}(\overline {X}({e_{{\vartheta}}})\in \dd z)
\,\pi(\dd x)}\\
&=&\frac{1}{{\vartheta}}\frac{\kappa({\vartheta},0)}{\kappa({\vartheta},\beta)}
\frac{\hat{\Phi}({\vartheta})}{\beta-\hat{\Phi}({\vartheta})}
\int_0^\infty\left( e^{-(\alpha+\hat{\psi}({\vartheta}))x}-e^{-(\alpha+\beta)x}\right)\,\pi(\dd x)\\
&=&\frac{\hat{\psi}'(0+)}{{\vartheta}}\frac{\kappa({\vartheta},0)}{\kappa({\vartheta},\beta)}\;
\frac{\hat{\Phi}({\vartheta})}{\beta-\hat{\Phi}({\vartheta})}
\left(\frac{\alpha+\hat{\Phi}({\vartheta})}
{\hat{\psi}(\alpha+\hat{\Phi}({\vartheta}))}-\frac{\alpha+\beta}{\hat{\psi}(\alpha+\beta)}\right).
\end{eqnarray*}
The Wiener-Hopf factorization \cite[Section 6.5.2] {kyprianoubook}
and (\ref{kappy}) complete the proof. \halmos

\vb

A similar result can be derived for the spectrally negative case.

\begin{proposition}\label{p.joint.transformn} If $X$ is a spectrally negative L\'{e}vy process, then
\begin{equation}
L({\vartheta};\alpha,\beta)= \frac{\Phi({\vartheta})-\alpha-\Phi(0)}{\Phi({\vartheta})+\beta}\frac{\Phi(0)}{\alpha+\beta+\Phi(0)}\frac{1}{{\vartheta}-\psi(\alpha+\Phi(0))}.\label{eq.Lsp.neg}
\end{equation}
\end{proposition}

\proof Recall the well-kown fact that  $\pi(\dd
x)=\Phi(0)e^{-\Phi(0)x}\,\dd x$ by (\ref{LTtau}). Applying Thm.\
\ref{lemma:basic} and interchanging the order of integration,
\begin{eqnarray*}
L({\vartheta};\alpha,\beta)&=&\frac{\Phi(0)}{{\vartheta}}\frac{\kappa({\vartheta},0)}{\kappa({\vartheta},\beta)}
\int_0^\infty e^{-(\alpha +\beta +\Phi(0)) x}\left(\int_0^x e^{\beta z}
\hat{\PP}(\overline {X}({e_{{\vartheta}}})\in \dd z)\right)\dd x\\
&=&\,\frac{1}{{\vartheta}}\frac{\kappa({\vartheta},0)}{\kappa({\vartheta},\beta)}\frac{\Phi(0)}{\alpha+\beta+\Phi(0)}
\hat{\EE}e^{-(\alpha +\Phi(0))\overline {X}({e_{{\vartheta}}}) }\\
&=&\,\frac{1}{{\vartheta}}\frac{\kappa({\vartheta},0)}{\kappa({\vartheta},\beta)}\frac{\Phi(0)}{\alpha+\beta+\Phi(0)}
\EE e^{-(\alpha +\Phi(0))\underline {X}({e_{{\vartheta}}}) }.
\end{eqnarray*}
This gives by Eqns.\ (\ref{kappy}) and (\ref{supdual}), in
conjunction with the fact that
${\vartheta}=\kappa({\vartheta},0)\hat{\kappa}({\vartheta},0)$,
\begin{eqnarray*}
L({\vartheta};\alpha,\beta)
&=&\frac{1}{{\vartheta}}\frac{\kappa({\vartheta},0)}{\kappa({\vartheta},\beta)}\frac{\Phi(0)}{\alpha+\beta+\Phi(0)}
\frac{\hat{\kappa}(0,0)}{\hat{\kappa}({\vartheta}, \alpha+\Phi(0))}\\
&=&\,\frac{1}{{\vartheta}}\frac{\kappa({\vartheta},0)}{\kappa({\vartheta},\beta)}\frac{\Phi(0)}{\alpha+\beta+\Phi(0)}
\frac{\hat{\kappa}({\vartheta},0)}{\hat{\kappa}({\vartheta}, \alpha+\Phi(0))}\\
&=&\,\frac{\Phi({\vartheta})-\alpha-\Phi(0)}{\Phi({\vartheta})+\beta}\frac{\Phi(0)}{\alpha+\beta+\Phi(0)}\frac{1}{{\vartheta}-\psi(\alpha+\Phi(0))}.
\end{eqnarray*}
This completes the proof.
\halmos

\section{Quasi-stationary distribution}\label{qsd}
In this section we use the Laplace transforms given
in (\ref{eq.Lsp.pos}) and (\ref{eq.Lsp.neg})  to identify the  quasi-stationary distribution $\mu(\dd x,\dd y)$
for the spectrally one-sided cases.
\subsection{Spectrally positive L\'evy process}\label{sec:pos}

We impose the following additional assumptions:
\newline {\bf [SP1]} \:\: There exists ${{\vartheta}_+}<0$ such that
\begin{itemize}
\item $\hat\psi({\vartheta})< \infty$  for ${{\vartheta}_+}<{\vartheta}$,
\item    $\hat{\psi}({\vartheta})$ attains its strictly negative minimum at ${\vartheta}{\s}<0$, where
${{\vartheta}_+}<{\vartheta}{\s}<0$ (and hence
$\hat{\psi}'({\vartheta}{\s})=0$).
\end{itemize}
Denote $\zeta{\s}:=\hat{\psi}({\vartheta}{\s})<0.$ Note that the
function $\hat\Phi$ can be considered in the complex domain. It is
clearly analytic for $\Re({\vartheta})>\zeta{\s}$. But, as it will
turn out, more is required to obtain the quasi-stationary
distribution.
\newline {\bf [SP2]} \:\: One can extend analytically $L({\vartheta};\alpha,\beta)$ into ${\mathscr G}_{\zeta{\s}}(\psi)$ for some $\pi/2<\psi\le
\pi$.

%$L({\vartheta};\alpha,\beta)$ fulfills
%condition (A1) and (A2) of Theorem \ref{t.tauberian}.

\begin{example} \rm
Since $\hat\Phi(\vartheta)$ is the Laplace exponent of a subordinator (viz.\
a first passage time process), we have the following spectral
representation:
\begin{equation}\label{reprhatphi}\hat\Phi({\vartheta})=d_+{\vartheta}+
\int_0^\infty (1-e^{-{\vartheta} x})\,\Pi_+(\dd x),\end{equation}
and $\int_0^\infty (x\land 1)\Pi_+(\dd x)<\infty$ \cite[Exercise
2.11]{kyprianoubook}. From its definition we see that $\zeta{\s}$
must be a singular point of $\hat\Phi({\vartheta})$. Moreover, if
there exists a density of $\Pi_+$ which is of semixponential type, then
from  Prop.~\ref{a1semexp} it follows then that $\hat{\Phi}$
is analytic in ${\mathscr G}_{\zeta{\s}}(\phi) $ and assumption
[SP2] is satisfied. In particular, assumption [SP2] is for example
satisfied for
 \[\Pi_+(\dd x)=e^{\zeta{\s} x}x^{\alpha}\,\dd x\,{\bf1}_{\{x>0\}},\] for $\alpha>-2$.
Clearly, then $\gamma({\vartheta})=\zeta{\s}
\cos{\vartheta}$.\hfill$\diamondsuit$
\end{example}

\begin{theorem}\label{specposmain}
If $X$ is a spectrally positive L\'evy process satisfying conditions {\em [SP1-SP2]}, then
\[\mu(\dd x,\dd y)=Q_+
\left(-\hat\psi'(0+)\frac{\hat\psi(\vt)-\vt\hat\psi'(\vt)}{(\hat\psi(\vt))^2}\right)
 e^{{\vartheta}{\s} (y-x)}x\,V_{-{\vartheta}{\s}}(y)\,\pi(\dd x)\dd y\,\;
{\boldsymbol 1}_{\{x\ge 0,y\ge 0\}},\] where $Q_+:=(\int_0^\infty
e^{{\vartheta}{\s} z}V_{-{\vartheta}{\s}}(z)\,{\rm
d}z)^{-1}$.\end{theorem}

\begin{corollary}
We have
\[
\mu^{\rm QS}_L (\dd x)=
\left(-\hat\psi'(0+)\frac{\hat\psi(\vt)-\vt\hat\psi'(\vt)}{(\hat\psi(\vt))^2}\right){e^{-{\vartheta}{\s}
x}x\,\pi(\dd x){\bs 1}_{\{x\geq 0\}}}\] and
\[\mu^{\rm QS}_R(\dd y)= Q_+{e^{{\vartheta}{\s} y}V_{-{\vartheta}{\s}}(y)\,\dd y
{\bs 1}_{\{y\geq 0\}}}.\]\end{corollary}

Before we prove Thm.\ \ref{specposmain}, we first present a few
facts. Let $k{\s}:=\sqrt{{2}/{ \hat{\psi}^{''}({\vartheta}{\s})}}.
$

\begin{lemma}\label{l.expansion1} Under {\em [SP1-SP2]},
$$\hat\Phi({\vartheta})={\vartheta}{\s}+ k{\s}({\vartheta}-\zeta^{*})^{1/2}+o(({\vartheta}-\zeta^{*})^{1/2})$$
as ${\vartheta}\downarrow\zeta{\s}$.
\end{lemma}

\proof From a Taylor series expansion and the condition that
$\hat\psi'({\vartheta}{\s})=0$, we have \[
\hat{\psi}({\vartheta})-\hat{\psi}({\vartheta}{\s})=
\frac{({\vartheta}-{\vartheta}{\s})^2}{2}
\psi''({\vartheta}{\s})+o(({\vartheta}-{\vartheta}{\s})^2).
\]
After some rearranging, it is obtained that
\[{{\vartheta}-{\vartheta}{\s}}\\
=\sqrt{\frac{2}{{\hat\psi}''({\vartheta}{\s})}}
\sqrt{\hat\psi({\vartheta})-\hat\psi({\vartheta}{\s})}+o(({\vartheta}-{\vartheta}{\s})^2).\]
We now substitute  ${\vartheta}=\hat\Phi(s)$, and use
$\hat{\psi}(\hat{\Phi}(s))=s$ to complete the proof. \halmos

\begin{proposition}\label{LTqspos}
Under {\em [SP1-SP2]} we have
\begin{eqnarray}
\lefteqn{\lim_{t\to\infty}\EE_{\pi}[e^{-\alpha Q(0)-\beta Q(t)}\,|\,T>t]}\nonumber\\
&=& \left(\zeta{\s}\cdot
\frac{\hat{\psi}(\alpha+{\vartheta}{\s})-(\alpha+{\vartheta}{\s})\hat{\psi}'(\alpha+{\vartheta}{\s})}
{\hat{\psi}^2(\alpha+{\vartheta}{\s})}
\right)\left(\frac{\zeta{\s}}{\zeta{\s}-\hat{\psi}(\beta)}\right)\;.\label{p.qs.formula}
\end{eqnarray}
\end{proposition}

\proof By [SP1-SP2] and Prop.\ \ref{p.joint.transform} $L({\vartheta}, \alpha,
\beta)$, as given in (\ref{eq.Lsp.pos}) as a function of
${\vartheta}$, is analytic in $\mathscr{G}_{\zeta{\s}}(\phi)$ for
$\pi/2<\phi\leq \pi$ when
$\hat{\psi}(\alpha+\hat{\Phi}(z))$ is
analytical there.
Recall that $\hat{\Phi}$ is analytic in this
region and note that $\hat{\psi}(\alpha+\hat{\Phi}(z))$
is analytical there since $\hat{\psi}(\hat{\Phi}(z))=z$ is
analytical in this region.
Thus condition (A1) of Thm.\ \ref{t.tauberian} is satisfied.
\begin{comment}
As long as in some neighborhood of $\hat{\Phi}(z_0)$ we
have $\hat{\Phi}^\prime(z)\neq 0$, the function $\hat{\psi}$ is
also analytic there. Thus, by the monodromy rule $\hat{\psi}$ is
analytical in the set of values of $\hat{\Phi}$ in
$\mathscr{G}_{\zeta{\s}}(\phi)$ when $\hat{\Phi}^\prime$ does not
disappear there. Now, this derivative could disappear only in the
region where $\hat{\psi}(\hat{\Phi}(z))$ is not analytical, that
is around $0$ where we should have $\hat{\psi}(\hat{\Phi}(z))=z$.
In fact for this region $\hat{\Phi}(z)$ converges to zero, and
$\hat{\psi}(\alpha+\hat{\Phi}(z))$ is still analytical as long as
for real $\alpha$ we have $\alpha >0$. \end{comment}

To check that condition (A2) of Thm.\ \ref{t.tauberian} holds for
$L({\vartheta}, \alpha, \beta)$, it suffices to prove that
\begin{equation}\label{mod}
\left|\frac{1}{z}\frac{\alpha+\hat{\Phi}(z)}{\hat{\psi}(\alpha+\hat{\Phi}(z))}\right|\end{equation}
tends to $0$ for $z\in \mathcal{G}_{\zeta{\s}}(\phi)$ tending to
$\infty$ two-dimensionally, that is in particular for
$\Im\,z\to\pm\infty$. From (\ref{reprhatphi}) it follows that
$b:=\Im\,\hat{\Phi}(z)\to\pm\infty$. Similarly one can prove that
either $a:=\alpha+\Re\,\hat{\Phi}(z)\leq 0$ is constant or tends
to $-\infty$. Then
\begin{eqnarray*}&&\Im\,\frac{\hat{\psi}(a+b\imi)}{|a+b\imi |}=\frac{1}{\sqrt{a^2+b^2}}\left(\hat{\psi}^\prime(0)b+ab\sigma^2+\int_\infty^0(\sin bx-b1_{\{|x|\leq 1\}}x)\hat{\Pi}(\dd x)
\right)\\&&\qquad\sim \hat{\psi}^\prime(0)-
\int_{-1}^0x\hat{\Pi}(\dd x)+a \sigma^2,
\end{eqnarray*}
which is either bounded or tends to $-\infty$ when $a$ is constant
or $a\to-\infty$. Taking into account the term $1/z$ in
(\ref{mod}), this completes the verification of condition (A2).
We will check now that also condition (A3) of Thm.\ \ref{t.tauberian} is
satisfied.
Now using  Lemma
\ref{l.expansion1}
we write
\begin{eqnarray*}
\hat\psi(\alpha+\hat{\Phi}({\vartheta}))&=&
\hat{\psi}(\alpha+{\vartheta}{\s}+k{\s}({\vartheta}-\zeta{\s})^{1/2}+o(({\vartheta}-\zeta{\s})^{1/2})\\&=&
\hat{\psi}(\alpha+{\vartheta}{\s})+\hat{\psi}^{'}(\alpha+{\vartheta}{\s})k{\s}({\vartheta}-\zeta{\s})^{1/2}+o(({\vartheta}-\zeta{\s})^{1/2}).
\end{eqnarray*}
Hence by Prop.\ \ref{p.joint.transform}  (for
${\vartheta}\downarrow \zeta{\s}$) we have, for some $\bar K$,
\begin{eqnarray*}
\lefteqn{L({\vartheta};\alpha,\beta)=
\frac{\hat{\psi}'(0+)}{{\vartheta}-\hat{\psi}(\beta)}\left(\frac{\alpha+\beta}{\hat{\psi}(\alpha+\beta)}
-
\frac{\alpha+{\vartheta}{\s}+k{\s}({\vartheta}-\zeta{\s})^{1/2}+o(({\vartheta}-\zeta{\s})^{1/2})}
{\hat{\psi}(\alpha+{\vartheta}{\s}+\gamma({\vartheta}-\zeta{\s})^{1/2}+o(({\vartheta}-\zeta{\s})^{1/2})}\right)}\\
&=&\bar K-\frac{\hat{\psi}'(0+)}{{\vartheta}-\hat{\psi}(\beta)}\;
\frac{\alpha+{\vartheta}{\s}+k{\s}({\vartheta}-\zeta{\s})^{1/2}+o(({\vartheta}-\zeta{\s})^{1/2}}
{\hat{\psi}(\alpha+{\vartheta}{\s})+\hat{\psi}^{'}(\alpha+{\vartheta}{\s})k{\s}({\vartheta}-\zeta{\s})^{1/2}+o(({\vartheta}-\zeta{\s})^{1/2})}\\
&=&\bar K-\frac{\hat{\psi}'(0+)}{{\vartheta}-\hat{\psi}(\beta)}\;
\frac{\alpha+{\vartheta}{\s}+k{\s}({\vartheta}-\zeta{\s})^{1/2}+o(({\vartheta}-\zeta{\s})^{1/2}}
{\hat{\psi}(\alpha+{\vartheta}{\s})+\hat{\psi}^{'}(\alpha+{\vartheta}{\s})k{\s}({\vartheta}-\zeta{\s})^{1/2}+o(({\vartheta}-\zeta{\s})^{1/2})}
\\
&&\hspace{2cm}\times \frac
{\hat{\psi}(\alpha+{\vartheta}{\s})-\hat{\psi}^{'}(\alpha+{\vartheta}{\s})k{\s}({\vartheta}-\zeta{\s})^{1/2}+o(({\vartheta}-\zeta{\s})^{1/2})}
{\hat{\psi}(\alpha+{\vartheta}{\s})-\hat{\psi}^{'}(\alpha+{\vartheta}{\s})k{\s}({\vartheta}-\zeta{\s})^{1/2}+o(({\vartheta}-\zeta{\s})^{1/2})}\;.
\end{eqnarray*}
Thus we obtain that, for some $K$,
\[L({\vartheta};\alpha,\beta)
= K-\frac{\hat{\psi}'(0+)k{\s}}{\zeta{\s}-\hat{\psi}(\beta)}
\frac{\hat{\psi}(\alpha+{\vartheta}{\s})-(\alpha+{\vartheta}{\s})\hat{\psi}'(\alpha+{\vartheta}{\s})}
{\hat{\psi}^2(\alpha+{\vartheta}{\s})}
({\vartheta}-\zeta{\s})^{1/2}+o(({\vartheta}-\zeta{\s})^{1/2}).\]
Conclude by invoking `Heaviside' that
\[\EE_\pi[e^{-\alpha Q(0)-\beta Q(t)},T>t]=
\frac{\hat{\psi}'(0+)k{\s}}{\zeta{\s}-\hat{\psi}(\beta)}
\frac{\hat{\psi}(\alpha+{\vartheta}{\s})-(\alpha+{\vartheta}{\s})\hat{\psi}'(\alpha+{\vartheta}{\s})}
{\hat{\psi}^2(\alpha+{\vartheta}{\s})}\frac{t^{-3/2}}{\Gamma(-1/2)}e^{\zeta{\s}t}(1+o(1)).\]
By setting $\alpha=\beta=0$ we have
\[\PP(T>t)=\frac{\hat{\psi}'(0+)k{\s}}{(\zeta{\s})^2}
\frac{t^{-3/2}}{\Gamma(-1/2)}e^{\zeta{\s} t}(1+o(1)).\] It is now
seen that Eqn.\  (\ref{p.qs.formula}) holds and the proof is
completed. \halmos

\vb

{\it Proof of Theorem \ref{specposmain}.}
From Prop.\ \ref{LTqspos} it follows that
 \[\tilde{\mu}(\alpha,\beta)=\int_0^\infty\int_0^\infty e^{-\alpha x}e^{-\beta y}\,
 \mu(\dd x,\dd y)=\tilde{A}_+(\alpha)\tilde{B}_+(\beta),\]
 where
 \[\tilde A_+(\alpha):={\zeta\s}\cdot\frac{\hat{\psi}(\alpha+{\vartheta}{\s})-(\alpha+{\vartheta}{\s})\hat{\psi}'(\alpha+{\vartheta}{\s})}
{\hat{\psi}^2(\alpha+{\vartheta}{\s})},\:\:\:\:\:\tilde
B_+(\beta):=\frac{\zeta\s}{\zeta{\s}-\hat\psi(\beta)}.\]

\noindent $\bullet$\:\: By the Pollaczek-Khintchine formula
(\ref{pollaczek}) and \cite[(4.14), p. 101]{kyprianoubook} applied
for the dual, we derive
\[\int_0^\infty e^{-(\alpha+\vartheta{\s}) x}x\;\pi(\dd x)=
\left.
-\tilde{\pi}^\prime({\vartheta})\right|_{{\vartheta}=\alpha+{\vartheta}{\s}}=
-\hat{\psi}^\prime(0+) \frac{ \hat{\psi}(\alpha+{\vartheta}{\s})-
(\alpha+{\vartheta}{\s})\hat{\psi}'(\alpha+{\vartheta}{\s})}
{\hat{\psi}^2(\alpha+{\vartheta}{\s})}.\] Hence, $\mu^{\rm
QS}_L(\cdot)$ has the desired form.

$\bullet$\:\: The dual process $\hat{X}$ is spectrally negative,
so that $\hat{V}_{-{\vartheta}{\s}}(y)=y$ and
$\hat{\kappa}_{-{\vartheta}{\s}}(0,{\vartheta})={\vartheta}$ for
all ${\vartheta}\geq 0$. The Wiener-Hopf factorization gives (up
to a multiplicative constant $k$ that relates to the normalization
of the local time) that under $\PP^{-{\vartheta}{\s}}$ for all
${\vartheta}\in \mathbb{R}$  we have
\begin{equation}
\label{WHtheta0} \psi_{-{\vartheta}{\s}}({\vartheta})
=-k{\vartheta} \kappa_{-{\vartheta}{\s}}(0, -{\vartheta})
\end{equation}
for all ${\vartheta}\leq -{{\vartheta}_+}$. From (\ref{newpsi})
and (\ref{LTV}) we have that
\begin{eqnarray*}
\lefteqn{\int_{0}^\infty e^{-\beta y} e^{{\vartheta}{\s}
y}{V}_{-{\vartheta}{\s}}(y)\;\dd y =
\frac{1}{(\beta-{\vartheta}{\s})\kappa_{-{\vartheta}{\s}}(0,\beta-{\vartheta}{\s})}}\nonumber\\&&
=\frac{k}{\psi_{-{\vartheta}{\s}} (-\beta + {\vartheta}{\s})}
=\frac{k}{\psi(-\beta)-\psi(-{\vartheta}{\s})} = \frac{k}{
\hat{\psi}(\beta)- \zeta{\s}}.
\end{eqnarray*}
Conclude that $\mu^{\rm QS}_R(\cdot)$ has the desired form, which
completes the proof. \halmos.

\begin{remark}\label{r1}\rm
The transform $\tilde{A}_+(\cdot)$  can be used to interpret the
quasi-stationary distributions.   Because of Pollaczek-Khinchine,
\[\frac{\alpha+{\vartheta}\s}{\hat\psi(\alpha+{\vartheta}\s)}\frac{\hat\phi({\vartheta}\s)}{{\vartheta}\s}\]
is a Laplace transform (i.e., corresponding to an exponentially
twisted version of the steady-state workload). In addition, by
virtue of \cite[Lemma 3.5]{abdel},
\[\frac{2}{\hat\psi'(0)}\frac{\alpha\hat\psi'(\alpha)-\hat\psi(\alpha)}{\alpha^2}\]
is a Laplace transform, and therefore also its
$\vartheta\s$-twisted version
\begin{equation}\label{LT3}\frac{(\alpha+{\vartheta}\s)\hat\psi'(\alpha+{\vartheta}\s)
-\hat\psi(\alpha+{\vartheta}\s)}{(\alpha+{\vartheta}\s)^2}
\frac{({\vartheta}\s)^2}{{\vartheta}\s\hat\psi'({\vartheta}\s)-\hat\psi({\vartheta}\s)}.\end{equation}
This reasoning indicates that, conditional on a long busy period,
$Q(0)$ is distributed as the sum of three independent random
variables. Two of these are distributed as the ${\vt}\s$-twisted
version of the steady-state workload, while a third has transform
(\ref{LT3}).\hfill$\diamondsuit$
\end{remark}

\subsection{Spectrally negative L\'evy process}\label{sec:neg}
Like for the spectrally positive case, also in the spectrally
negative case we need to impose additional assumptions to find the
quasi-stationary distribution.
\newline {\bf [SN1]} \:\: There exists ${{\vartheta}_-}>0$ such that
\begin{itemize}
\item $\psi({\vartheta})< \infty$  for $0<{\vartheta}<{{\vartheta}_-}$,
\item    $\psi({\vartheta})$ attains its strictly negative minimum at ${\vartheta}{\s}>0$, where
$0<{\vartheta}{\s}<{{\vartheta}_-}$ (and hence
$\psi'({\vartheta}{\s})=0$).
\end{itemize}
{\bf [SN2]} \:\: $\Phi$ is analytical in
$\mathscr{G}_{\zeta{\s}}(\phi)$ for $ \pi/2< \phi\leq\pi$, where
$\zeta{\s}:=\psi({\vartheta}{\s})<0.$

\begin{example} \rm
Since for spectrally negative L\'evy process $\Phi(\vartheta)$ is the Laplace exponent of a subordinator (viz.\
a first passage time process), the spectral
representation
\begin{equation}\label{reprhatphi2}\Phi({\vartheta})=d_-{\vartheta}+
\int_0^\infty (1-e^{-{\vartheta} x})\,\Pi_-(\dd x),\end{equation}
applies, with
$\int_0^\infty (x\land 1)\Pi_-(\dd x)<\infty$; cf.\ (\ref{reprhatphi}). This means that if
there exists a density of $\Pi_-$ which is of semixponential type, then
Prop.~\ref{a1semexp} entails that $\Phi$
is analytic in ${\mathscr G}_{\zeta{\s}}(\phi) $ and hence assumption
[SN2] is satisfied. \hfill$\diamondsuit$
\end{example}

\begin{theorem}\label{specnegmain}
If $X$ is a spectrally negative L\'evy process satisfying conditions {\em [SN1-SN2]}, then
\begin{eqnarray*}
\mu(\dd x,\dd y)=Q_-({\vartheta}{\s})^2ye^{-{\vartheta}{\s} (x+y)}
e^{-\Phi(0)x}\hat{V}_{{\vartheta}{\s}}(x)\,\dd x\,\dd y\; {\boldsymbol
1}_{\{x\geq 0,y\geq 0\}},
\end{eqnarray*}
where
%$D=(\zeta{\s})^2/\hat{\psi}(0+)$.
$Q_-:=(\int_0^\infty
e^{-(\Phi(0)+{\vartheta}{\s})z}V_{{\vartheta}{\s}}(z)\,\dd
z)^{-1}$.
\end{theorem}

\begin{corollary}
We have
\begin{eqnarray*}
\mu^{\rm QS}_L(\dd y)=
Q_-e^{-(\Phi(0)+{\vartheta}{\s})x}\hat{V}_{{\vartheta}{\s}}(x)\,\dd
x{\boldsymbol 1}_{\{x\geq 0\}}
\end{eqnarray*}
and
\[\mu^{\rm QS}_R (\dd y)= ({\vartheta}{\s})^2ye^{-{\vartheta}{\s} y}
\dd y{\boldsymbol 1}_{\{y\geq 0\}}\;.
\]
\end{corollary}
Observe that $\mu^{\rm QS}_R (\cdot)$ corresponds with an
Erlang(2) distribution. The proof of these results is based on the
following lemma, which is proven as  Lemma \ref{l.expansion1}.

\begin{lemma}\label{l.expansion2} Under {\em [SN1-SN2]},
\[\Phi({\vartheta})={\vartheta}{\s}+ k{\s}({\vartheta}-\zeta^{*})^{1/2}+o(({\vartheta}-\zeta^{*})^{1/2})\]
as ${\vartheta}\downarrow\zeta{\s}$, where $k{\s}:=\sqrt{{2}/{
\psi^{''}({\vartheta}{\s})}}.$
\end{lemma}

{\it Proof of Theorem \ref{specnegmain}.} Note that all
assumptions of Thm.\ \ref{t.tauberian} are satisfied  by Prop.\
\ref{p.joint.transformn}.
In particular, as $|z|$ tends to infinity in $\mathscr{G}_{\zeta{\s}}(\phi)$ function
$|\Phi(z)|$ is either bounded or tends to infinity. In both cases condition (A2) is satisfied.
Moreover, for some $K$,
\[L(\alpha,\beta;\vartheta)=K-\frac{\Phi(0)k{\s}}{({\vartheta}{\s}+\beta)^2}\frac{1}{\psi(\alpha+\Phi(0))-\zeta{\s}}
({\vartheta}-\zeta{\s})^{1/2}+o(({\vartheta}-\zeta{\s})^{1/2})\]as
${\vartheta}\downarrow \zeta{\s}$, and we can conclude by
`Heaviside'  that
\[\EE_\pi[e^{-\alpha Q(0)-\beta Q(t)},T>t]=\frac{\Phi(0)k{\s}}{({\vartheta}{\s}+\beta)^2}\frac{1}{\psi(\alpha+\Phi(0))-\zeta{\s}}
\frac{t^{-3/2}}{\Gamma(-1/2)}e^{-\zeta{\s}t}(1+o(1))\]as
$t\to\infty$. Therefore,\[
\tilde{\mu}(\alpha,\beta)=\lim_{t\to\infty}\EE_{\pi}[e^{-\alpha
Q(0)-\beta Q(t)}\,|\,T>t]=
\tilde{A}_-(\alpha)\tilde{B}_-(\beta),\] where
\[\tilde{A}_-(\alpha):=\frac{-\zeta\s}{\psi(\alpha+\Phi(0))-\zeta{\s}},\qquad
\tilde{B}_-(\beta):=\frac{({\vartheta}{\s})^2}{({\vartheta}{\s}+\beta)^2}\;.\]
It is not hard to see that the proposed density indeed corresponds
with this transform. \halmos

\section{Examples}
In this section we illustrate our theory by means of a number of examples.
We indicate for what L\'evy processes our assumptions are fulfilled, and for a few of those processes
we perform the computations.

According to Vigon's theory of {\it philanthropy}  \cite{23},  a (killed)
subordinator is called a philanthropist if its L\'evy measure has a
decreasing density on ${\mathbb R}_+$. Moreover, given any two
subordinators $H_1$ and $H_2$ which are philanthropists, providing
that at least one of them is not killed, there exists a L\'evy
process $X$ such that $H_1$ and $H_2$ have the same law as the
ascending and descending ladder height processes of $X$,
respectively. Suppose we denote the killing rate, drift
coefficient and L\'{e}vy measures of $H_1$ and $H_2$ by the
respective triples $(b, \delta,\Pi_{H_1})$ and $(\hat{b},
\hat{\delta},\Pi_{H_2})$. Then \cite{23} shows that the L\'{e}vy
measure of $X$ satisfies the following identity
\[
\Pi(x,\infty)=\int_0^\infty \Pi_{H_2}(u,\infty)\,\Pi_{H_1}(x+{\rm d}u)+\hat{\delta}\pi_{H_1}(x)+\hat{b}\Pi_{H_1}(x,\infty),\qquad x>0,
\]
where $\pi_{H_1}(x)$ is the density corresponding to $\Pi_{H_1}$. By symmetry, an obvious analogue of the above equation holds for the
negative tail $\Pi(-\infty,x)$, with $x<0$.

Choosing then e.g. $H_2(t)=t$ and $H_1(t)=\tau(t)$ with Laplace exponent $\hat{\Phi}$ and jump measure $\Pi_+$
being semiexponential, then, using the above construction, we can easily give examples of
a spectrally positive L\'{e}vy process $X$ satisfying conditions
[SP1-SP2]. Similarly, using the above method we can construct spectrally negative L\'{e}vy processes satisfying [SN1-SN2].

\vb

Usually these conditions can be verified in a straightforward manner, as we did in the examples below.

\begin{comment}
Check $\beta$-subordinators of \cite{sim18}. inverse gaussian etc.
To be completed... \footnote{\tt ZP: For a spectrally positive
L\'{e}vy process renewal function $V_{-{\vartheta}{\s}}(x)$ equals
$\hat{W}_{-{\vartheta}{\s}}(x) =e^{-\hat{\Phi}(\zeta{\s})x}
\hat{W}^{(\zeta{\s})}(x)$, hence is known as long we can identify
the scale function of the dual. Few further examples could be
given using Andreas papers...} \footnote{\tt MM: This part needs
to be rewritten...}
\end{comment}

\newcommand{\vr}{\varrho}

\begin{example}\rm
M/M/1 {\it  queue.} In this case
\begin{eqnarray}\label{CL}
 X(t) = \sum_{i=1}^{N(t)} \sigma_i-t,
\end{eqnarray}
where $\sigma_i$ (where $i=1,2,...$) are i.i.d.\ service times
that have an exponential distribution with mean $1/\nu$. The
arrival process is a homogeneous Poisson process $N(t)$ with rate
$\lambda$; it is assumed that $\varrho:=\lambda/\nu<1.$ We apply
the theory of Section \ref{sec:pos}.

\begin{comment}
Let us first concentrate on $\mu^{\rm QS}_L(\cdot)$. Then
\end{comment}
We have
\[\hat\psi(\eta)=\eta-\lambda\left(1-\frac{\nu}{\eta + \nu}\right)=\eta-
\frac{\lambda\eta}{\eta+\nu},\]
yielding ${\vt}\s=\sqrt{\lambda\nu}-\nu$, and
$\zeta\s=-(\sqrt{\nu}-\sqrt{\lambda})^2.$
Furthermore
$$\widehat{\Phi}(\eta)=\frac{{\eta} + \lambda -\nu  +\sqrt{({\eta} + \lambda-\nu)^2 + 4 {\theta}\nu}}{2 }$$
and hence assumptions [SP1-SP2] are satisfied.

Let us first concentrate on $\mu^{\rm QS}_L(\cdot)$. We  use Remark \ref{r1}.
Using that
\[{\eta\hat\psi'(\eta)-\hat\psi(\eta)}=\frac{\lambda{\eta^2}}{(\eta+\nu)^2},\]
we obtain
\[\frac{(\alpha+{\vartheta}\s)\hat\psi'(\alpha+{\eta}\s)
-\hat\psi(\alpha+{\vartheta}\s)}{(\alpha+{\vartheta}\s)^2}
\frac{({\vartheta}\s)^2}{{\vartheta}\s\hat\psi'({\vartheta}\s)-\hat\psi({\vartheta}\s)}
=\left(\frac{\sqrt{\lambda\nu}}{\alpha+\sqrt{\lambda\nu}}\right)^2.\]
This corresponds with the sum of two Exp($\sqrt{\lambda\nu}$)
random variables. Also,
\[\frac{\alpha+{\vartheta}\s}{\hat\psi(\alpha+{\vartheta}\s)}\frac{\hat\psi({\vartheta}\s)}
{{\vartheta}\s}=
(1-\sqrt{\vr})\frac{\alpha+\sqrt{\lambda\nu}}{\alpha+\sqrt{\lambda\nu}-\lambda}=
(1-\sqrt{\vr})\sum_{n=0}^\infty
(\sqrt{\vr})^n
\left(\frac{\sqrt{\lambda\nu}}{\alpha+\sqrt{\lambda\nu}}\right)^n,\]
corresponding with a shifted-Geom($\sqrt{\vr}$)-distributed number
of Exp($\sqrt{\lambda\nu}$) random variables. Conclude that
$\mu^{\rm QS}_L(\cdot)$ corresponds to the sum of $M$ independent
Exp($\sqrt{\lambda\nu}$) random variables, where
\[\PP(M=m) = (m-1)(\sqrt{\vr})^{m-2}(1-\sqrt{\vr})^2,\]
i.e., $M$ has a negative binomial distribution with parameters 2
and $\sqrt{\vr}$. A similar form is found for the general
light-tailed M/G/1 case.

Let us now study $\mu^{\rm QS}_R(\cdot)$. It is a matter of
straightforward calculus to find that
\[\tilde B_+(\beta) = (\nu+\beta)\left(\frac{\sqrt{\nu}-\sqrt{\lambda}}
{\beta+\sqrt{\nu}(\sqrt{\nu}-\sqrt{\lambda})}\right)^2.\] A
partial fraction expansion argument gives that this equals
\[(1-\sqrt{\vr})\frac{\nu-\sqrt{\lambda\nu}}{\nu-\sqrt{\lambda\nu}+\beta}+\sqrt{\vr}
\left(\frac{\nu-\sqrt{\lambda\nu}}{\nu-\sqrt{\lambda\nu}+\beta}\right)^2.\]
In other words, the quasi-stationary distribution at time $t$ (for
$t$ large) equals a mixture of an exponential and an Erlang(2)
distribution.
 \hfill$\diamondsuit$
\end{example}

\iffalse Then
$$\hat{\Phi}({\vartheta})=\frac{{{\vartheta}} + \lambda -\xi  +\sqrt{({{\vartheta}} + \lambda-\xi)^2 + 4 {{\vartheta}}\xi}}{2 }$$
and hence assumptions [SP1-SP2] are satisfied.
\footnote{My calculations:
$$\psi({\vartheta})=-{\vartheta}+\lambda(m({\vartheta})-1)$$
$$\hat{\psi}({\vartheta})={\vartheta}+\lambda(m(-{\vartheta})-1)$$
and for the exp. distr.
$$\hat{\psi}({\vartheta})={\vartheta}-\frac{\lambda{\vartheta}}{\xi+{\vartheta}}$$
$$\hat{\psi}^{'}({\vartheta})=1-\frac{\lambda\xi}{(\xi+{\vartheta})^2}$$
$${\vartheta}{\s}=-\xi+\sqrt{\lambda\xi}$$
$$({{\vartheta}} + \lambda-\xi)^2 + 4 {{\vartheta}}\xi=0$$ fpr
$${\vartheta}=-(\lambda^{1/2}\pm \xi^{1/2})^2$$
and so $\zeta{\s}=-(\lambda^{1/2}- \xi^{1/2})^2$ and
$\hat{\Phi}(z)$ is analytic in $\mathbb{C}-(-\infty,\zeta{\s}]$.
Now
$$({\vartheta}+\lambda-\xi)^2+4{\vartheta}\xi=0$$
has the biger solution $\zeta{\s}=-(\lambda^{1/2}-\xi^{1/2})$.}
Hence $\zeta{\s}=-(\lambda^{1/2}- \xi^{1/2})^2$ and
$\hat{\Phi}(z)$ is analytic in
$\mathbb{C}-(-\infty,\zeta{\s}]$.\fi

\begin{example}\rm
{\it Linear Brownian motion.} In this case $X(t)=\sigma
B(t)-t,$ where $\sigma>0$ and $B(t)$ is a standard Brownian
motion. Remark that this process is spectrally positive {\it and}
spectrally negative, so we can use both Thm.\ \ref{specposmain}
and Thm.\ \ref{specnegmain}.

Let us first see what the spectrally positive results would give.
It is not hard to check that
\[\hat{\psi}({\vartheta})={\vartheta}+\frac{\sigma^2{\vartheta}^2}{2},\]
so that, in the setting of Section \ref{sec:pos}, $\vartheta\s=
-1/\sigma^2$ and $\zeta\s=-1/(2\sigma^2).$ It is a matter of
straightforward computations now to obtain that
\[\tilde{A}_+(\alpha)=\left(\frac{1/\sigma^2}{1/\sigma^2 +\alpha}\right)^2,\:\:\:\:
\tilde{B}_+(\beta)=\left(\frac{1/\sigma^2}{1/\sigma^2
+\beta}\right)^2.\] Conclude that the quasi-stationary
distributions of $Q(0)$ and $Q(t)$ ($t$ large) are both Erlang(2)
with mean $2/\sigma^2$, whereas the stationary workload has an
exponential distribution with mean $1/(2\sigma^2).$ (In the
decomposition of Remark \ref{r1}, the first two random variables
have exponential distributions with mean $1/\sigma^2$, the third
is equal to 0). Interestingly, the relation with the Erlang(2)
distribution has also been observed in, e.g.,
\cite{Iglehart,Martmart, PR}.

The same result can be obtained by using the results from Section
\ref{sec:neg}. Now $\vt\s=1/\sigma^2$ and
$\zeta\s=-1/(2\sigma^2).$ It is easily checked that
$\Phi(0)=2/\sigma^2$. As expected, we obtain
$\tilde{A}_-(\alpha)=\tilde{A}_+(\alpha)$ and
$\tilde{B}_-(\beta)=\tilde{B}_+(\beta)$.

\newcommand{\NN}{_{\rm N}}

In fact, in this case the quasi-stationarity distributions can be
found in an explicit manner. For ease we restrict ourselves to
studying just $\mu^{\rm QS}_L(\cdot)$; we do so by investigating
the density
\[\frac{\rm d}{{\rm d} q}\PP_\pi(Q(0)\le q\mid T>t)=:f_t(q).\]
We rely on the standard equality
\[\PP_q(T>t)= \Psi\NN\left(\frac{t-q}{\sqrt{t}}\right)-e^{2q}\Phi\NN\left(\frac{-t-q}{\sqrt{t}}\right),\]
and the fact that $Q(0)$ (unconditioned) has an exponential
distribution with mean $\frac{1}{2}$; here, $\Phi\NN(\cdot)$
denotes the distribution function of a standard Normal random
variable, where $\Psi\NN(x):=1-\Phi\NN(x)$ is its tail. It is
known that, as $x\to\infty$,
\begin{equation}\label{assphi}
\Psi\NN(x)\sim
\left(\frac{1}{x}-\frac{1}{x^3}+\frac{3}{x^5}\right)\frac{{1}}{\sqrt{2\pi}}e^{-\frac{1}{2}x^2}.\end{equation}
Let us first determine $\PP_\pi(T>t)$ (where $Q(0)$ has an
exponential distribution with mean $\frac{1}{2}$), which can
evidently be rewritten as
\[1-\int_0^\infty 2e^{-2q}\Phi\NN\left(\frac{t-q}{\sqrt{t}}\right){\rm d}q-2\int_0^\infty \Phi\NN\left(\frac{-t-q}{\sqrt{t}}\right){\rm d}q.\]
Consider the first integral of the previous display. It can be
evaluated as
\begin{eqnarray*}
\lefteqn{\int_0^\infty\int_{-\infty}^{\frac{t-q}{\sqrt{t}}}2e^{-2q}\cdot \frac{1}{\sqrt{2\pi}}e^{-\frac{1}{2}y^2}{\rm d}y{\rm d}q=\int_{-\infty}^{\sqrt{t}}\int_0^{t-y\sqrt{t}}  2e^{-2q}\cdot \frac{1}{\sqrt{2\pi}}e^{-\frac{1}{2}y^2}{\rm d}q{\rm d}y}\\
&=&\int_{-\infty}^{\sqrt{t}}\left(1-e^{-2(t-y\sqrt{t})}\right)
\cdot \frac{1}{\sqrt{2\pi}}e^{-\frac{1}{2}y^2}{\rm
d}y=\Phi\NN(\sqrt{t})-\Phi\NN(-\sqrt{t}).
\end{eqnarray*}
Likewise, the second integral can be rewritten as
\begin{eqnarray*}
\lefteqn{\int_0^\infty\int_{-\infty}^{\frac{-t-q}{\sqrt{t}}}2\cdot \frac{1}{\sqrt{2\pi}}e^{-\frac{1}{2}y^2}{\rm d}y{\rm d}q=\int_{-\infty}^{-\sqrt{t}}\int_0^{-t-y\sqrt{t}}  2\cdot \frac{1}{\sqrt{2\pi}}e^{-\frac{1}{2}y^2}{\rm d}q{\rm d}y}\\
&=&\int_{-\infty}^{-\sqrt{t}}(-2t-2y\sqrt{t})  \cdot
\frac{1}{\sqrt{2\pi}}e^{-\frac{1}{2}y^2}{\rm
d}y=-2t\Phi\NN(-\sqrt{t})+2\sqrt{t}\cdot\frac{1}{\sqrt{2\pi}}e^{-\frac{1}{2}t}.
\end{eqnarray*}
We arrive at
\begin{equation}\label{den}\PP_\pi(T>t)=\Psi\NN(\sqrt{t})+\Phi\NN(-\sqrt{t})+2t\Phi\NN(-\sqrt{t})-2\sqrt{t}\cdot\frac{1}{\sqrt{2\pi}}e^{-\frac{1}{2}t}.\end{equation}
Using (\ref{assphi}) it is readily verified that, for $t$ large,
\[\PP_\pi(T>t)\sim\frac{4}{t\sqrt{t}}\cdot\frac{1}{\sqrt{2\pi}}e^{-\frac{1}{2}t}.\]
Also, it holds that
\begin{equation}\label{num}
\frac{\rm d}{{\rm d} q}\PP_\pi(Q_0\le q,
T>t)=2e^{-2q}\Psi\NN\left(\frac{t-q}{\sqrt{t}}\right)
-2\Phi\NN\left(\frac{-t-q}{\sqrt{t}}\right),\end{equation} so that we
now have an explicit expression for $f_t(q)$, viz.\ the ratio of
(\ref{num}) and (\ref{den}). Due to the asymptotic equivalence
(\ref{assphi}), Expression (\ref{num}) behaves for $t$ large as
\[\frac{2}{\sqrt{2\pi}}{\sqrt{t}}\left(\frac{1}{t-q}-\frac{1}{t+q}\right)\exp\left(-\frac{1}{2}\frac{(q+t)^2}{t}\right)
\sim
4\left(\frac{1}{t\sqrt{t}}\frac{1}{\sqrt{2\pi}}e^{-\frac{1}{2}t}\right)q
e^{-q}.\] We conclude that we again find that the quasi-stationary
distribution of $Q(0)$ is Erlang(2) with expected value $2$.
\hfill$\diamondsuit$ \iffalse Then
$\hat{\Phi}({\vartheta})=\frac{-1+\sqrt{1-2\sigma^2{\vartheta}}}{\sigma^2}$
and $\zeta{\s}=\frac{-1}{2\sigma^2}$. and also conditions
[SP1-SP2] are satisfied. Notice that ${\mathscr
G}_{\zeta{\s}}=\mathbb{C} -(-\infty,\zeta{\s}]$. Hence
$$\mu(dx,dy)=Dxye^{-\frac{x+y}{\sigma^2}}.$$
Now
$$D=1/(\int_0^\infty xe^{-x/\sigma^2}\,dx)^2=1/\sigma^8.$$
Mention that the Erlang(2) also comes up in, e.g., \cite{Martmart,
PR}. \footnote{My calculations:
$$\psi({\vartheta})=-{\vartheta}+\frac{\sigma^2{\vartheta}^2}{2}$$
so
$$\hat{\psi}({\vartheta})={\vartheta}+\frac{\sigma^2{\vartheta}^2}{2}.$$
Then
$$\hat{\psi}^{'}({\vartheta})=1+\sigma^2{\vartheta},\qquad \hat{\psi}^{''}({\vartheta})=\sigma^2$$
and
$\hat{\Phi}({\vartheta})=\frac{-1+\sqrt{1+2\sigma^2 {\vartheta}}}{\sigma^2}$.
Thus
$${\vartheta}^{*}=-1/\sigma^2,\qquad \zeta{\s}=-\frac{1}{\sigma^2}.$$
$$k^{*}=\left(\frac{2}{\sigma^2}\right)^{1/2}=\frac{\sqrt{2}}{\sigma}.$$
$$\pi(dx)=\frac{2}{\sigma^2}e^{-2y/\sigma^2}$$
  }\fi\end{example}

  \section*{Acknowledgment}
  The authors thank P.\ Glynn (Stanford University) for inspiring discussions.

\end{document}